\documentclass[preprint,11pt]{amsart}
\usepackage[margin=1in]{geometry}              % to include figures
\usepackage{amsmath}               % great math stuff
\usepackage{amsfonts}              % for blackboard bold, etc
\usepackage{amsthm}
\usepackage{enumerate}
\usepackage{esint}
\usepackage{xcolor}
\usepackage{soul}
\sethlcolor{yellow}

\numberwithin{equation}{section}

\newtheorem{thm}{Theorem}[section]
\newtheorem{defn}[thm]{Definition}
\newtheorem{lem}[thm]{Lemma}

\newtheorem{Rem}[thm]{Remark}

\def\XXint#1#2#3{{\setbox0=\hbox{$#1{#2#3}{\int}$}
     \vcenter{\hbox{$#2#3$}}\kern-.5\wd0}}

\begin{document}

\nocite{*}

\title[Local Boundedness]{Local Boundedness of a Direction Dependent Double Phase Nonlocal Elliptic Equation }


\author{Hamid EL Bahja}

\address{Hamid EL Bahja, The African Institute for Mathematical Sciences, Research and Innovation Centre, Rwanda}
\email{hamidsm88@gmail.com, helbahja@aims.ac.za}


\maketitle

\begin{abstract}
This paper investigates the local boundedness of weak solutions to a direction-dependent double-phase nonlocal elliptic equation. By employing refined energy estimates and De Giorgi-type techniques, we establish the local boundedness of these solutions.
\end{abstract}

\section{Introduction and statement of main results}
In this paper, we are concerned with regularity properties of the following anisotropic fractional nonlocal double phase equation of the form
\begin{equation}
   \mathcal{L}u=0~~~~~~~~\text{                    in}~~~~\Omega
\end{equation}
defined in a bounded domain $\Omega$ in $\mathbb{R}^{N}$ by
\begin{equation}
\begin{aligned}
\mathcal{L}u(x)
= \sum_{i=1}^{N} &\operatorname{P.V.}\!\int_{\mathbb{R}}
\Bigg[
\, s_{i}(1-s_{i})
\frac{|u(x)-u(x+he_{i})|^{p-2}\big(u(x)-u(x+he_{i})\big)}{|h|^{1+s_{i}p}}  \\
&+\, t_{i}(1-t_{i})\,a(x,x+he_{i})
\frac{|u(x)-u(x+he_{i})|^{q-2}\big(u(x)-u(x+he_{i})\big)}{|h|^{1+t_{i}q}}
\Bigg] \, dh 
\end{aligned}
\end{equation}
where 
\begin{equation}0<s_{i}\leq t_{i}<1<p\leq q<\infty\end{equation}
for all $i=1,\dots,N$, and $a:\mathbb{R}^{N}\times \mathbb{R}^{N}\rightarrow\mathbb{R}$ is a nonnegative, measurable, and bounded modulating coefficient such that
\begin{equation}
    0\leq a(x,y)\leq \|a\|_{L^{\infty}},~~x,y\in\mathbb{R}^{N}.
\end{equation}
The operator used in our study can be seen as a nonlocal analog of the following orthotopic double-phase operator:  
\[
\sum_{i=1}^{N} \frac{\partial}{\partial x_{i}} \left( \left| \frac{\partial u}{\partial x_{i}} \right|^{p-2} \frac{\partial u}{\partial x_{i}} \right) 
+ a(x) \sum_{i=1}^{N} \frac{\partial}{\partial x_{i}} \left( \left| \frac{\partial u}{\partial x_{i}} \right|^{q-2} \frac{\partial u}{\partial x_{i}} \right),
\]
as \( s, t \to 1 \), where \( a: \mathbb{R}^{N} \to \mathbb{R} \) is a negative modulating function. This class of operators, characterized by anisotropic and singular kernels of the form (1.2), naturally arises in various mathematical frameworks, including stochastic differential equations and potential theory \cite{bas,fri,kass1,kul}. Furthermore, they play a crucial role in the study of anisotropic stable jump‑Cox–Ingersoll–Ross (jump‑CIR) processes, which generalize classical jump-diffusion models by incorporating directional dependencies and nonlocal effects \cite{fri2}. The interplay between anisotropy, nonlocality, and double-phase growth in such models introduces significant analytical challenges, motivating further investigation into their qualitative and quantitative properties. It is worth mentioning that the regularity theory of the local anisotropic \( p \)-Laplacian has been studied in \cite{dib,ham1,ham2} and the references therein, providing a foundational understanding that informs the analysis of its nonlocal counterparts.

The regularity theory for nonlocal problems with fractional orders has been extensively studied over the last two decades. Caffarelli and Silvestre \cite{caf1} established the Harnack inequality for the fractional Laplacian \( (-\Delta)^s u = 0 \) using an extension technique. Later, Caffarelli, Chan, and Vasseur \cite{caf2} applied De Giorgi’s method to linear parabolic equations with general fractional kernels, proving the Hölder continuity of weak solutions. Extensive research has since been conducted on the regularity of nonlocal linear equations with fractional orders; see, for instance, \cite{ros,kass2,sil,fel} and references therein.

For nonlinear nonlocal problems, particularly those involving the fractional \( p \)-Laplacian, Di Castro, Kuusi, and Palatucci \cite{cas1,cas2} developed a nonlocal variant of De Giorgi’s approach, establishing local Hölder continuity and Harnack inequalities. Cozzi \cite{Coz} later extended these results to non-homogeneous equations by introducing fractional De Giorgi classes. These methods have been instrumental in addressing various problems, including obstacle problems \cite{kor} and measure data problems \cite{kuu}. For a broader discussion on these developments, we refer to \cite{ok,byu,din,nez} and the monographs \cite{buc,kuu2}.

Regarding isotropic nonlocal problems involving fractional double phase
operators, the corresponding regularity theory is still in a developing
stage and has been investigated mainly in stationary settings. In the
context of mixed local and nonlocal equations with double phase
structure, De Filippis and Mingione \cite{fil2} established gradient
regularity results, highlighting the variational nature of such models.
For purely nonlocal operators, De Filippis and Palatucci \cite{fil}
proved interior H\"older continuity for viscosity solutions of the
equation $Lu = f \in L^{\infty}$. Concerning weak solutions, Fang and
Zhang \cite{fan} obtained local H\"older continuity (without specifying
the exponent) for bounded weak solutions by means of a
De Giorgi--Nash--Moser type method under the condition $q s_2 \leq p s_1$
and a nonnegative bounded modulating coefficient, also showing that
these solutions are viscosity solutions. Subsequently, Byun et al.
\cite{byu2} extended these results to the case $p s_1 < q s_2$ by using a
logarithmic lemma combined with the De Giorgi--Nash--Moser iteration,
assuming that the modulating coefficient belongs to
$C^{0,\beta}(\mathbb{R}^N \times \mathbb{R}^N)$ with
$q s_2 - p s_1 < \beta < 1$.

Recent research has advanced the understanding of nonlocal equations with anisotropic and singular kernels. Chaker et al. \cite{cha1,cha2,cha3} investigated nonlocal operators featuring singular anisotropic kernels, focusing on anisotropic jump processes with direction-dependent stability indices. Their work established weak Harnack inequalities and Hölder continuity results for solutions, offering new insights into the interplay between anisotropy and singularity in nonlocal settings.

Let \(\vec{s} = (s_1, s_2, \dots, s_N)\) with \(s_i \in (0, 1)\) for \(i = 1, 2, \dots, N\), and let \(p > 1\) satisfy
\begin{equation}
p < p^{*}_{\bar{s}} = \frac{Np}{N - \bar{s}p}, \quad \text{where} \quad \bar{s} = \left( \frac{1}{N} \sum_{i=1}^N \frac{1}{s_i} \right)^{-1} \quad \text{and} \quad \bar{s}p < N.
\end{equation}
We introduce the homogeneous anisotropic Sobolev space $D^{\vec{s}, p}(\mathbb{R}^N, \mathbb{R})$. 
To give an unambiguous meaning to the pointwise differences that appear, we work with a Borel measurable representative of each equivalence class in $L^p(\mathbb{R}^N)$. 
The space is then defined as
\[
D^{\vec{s}, p}(\mathbb{R}^N, \mathbb{R})
=
\left\{
u \in L^p(\mathbb{R}^N)
\;:\;
\sum_{i=1}^N
\int_{\mathbb{R}^N}
\int_{\mathbb{R}}
s_i (1 - s_i)
\frac{|u(x) - u(x + h e_i)|^p}{|h|^{1 + s_i p}}
\, dh \, dx
< \infty
\right\}.
\]
For $u \in D^{\vec{s}, p}(\mathbb{R}^N,\mathbb{R})$, the function
\[
(x,h) \longmapsto 
s_i (1 - s_i)
\frac{|u(x) - u(x + h e_i)|^p}{|h|^{1 + s_i p}}
\]
belongs to $L^1(\mathbb{R}^N \times \mathbb{R})$ for each $i=1,\dots,N$. 
Consequently, by Fubini's theorem, the above double integral is well-defined and the order of integration may be interchanged. 
In particular, for each $i$, the inner $h$-integral\[
I_i(x):=\int_{\mathbb{R}}
s_i (1 - s_i)
\frac{|u(x) - u(x + h e_i)|^p}{|h|^{1 + s_i p}}
\, dh
\]
is finite for almost every $x\in\mathbb{R}^N$, and the map $x\mapsto I_i(x)$ is measurable. 
Accordingly, all pointwise expressions involving $u(x+he_i)$ in this paper are to be interpreted in the almost-everywhere sense.

We endow $D^{\vec{s}, p}(\mathbb{R}^N,\mathbb{R})$ with the norm
\[
\|u\|_{D^{\vec{s}, p}}
=
\|u\|_{L^p(\mathbb{R}^N)} 
+ 
\sum_{i=1}^N 
\left(
\int_{\mathbb{R}^N}
\int_{\mathbb{R}} 
s_i(1-s_i)
\frac{|u(x)-u(x+he_i)|^p}{|h|^{1+s_ip}}
\,dh\,dx
\right)^{1/p}.
\]

It follows from [\cite{cha3}, Theorem 1.1] that $D^{\vec{s}, p}(\mathbb{R}^N,\mathbb{R})$ is a normed space.

To address the anisotropy of the problem, we consider rectangles tailored to the fractional parameters \(s_1, s_2, \dots, s_N\). For \(\rho > 0\) and \(x \in \mathbb{R}^N\), define the rectangle
\begin{equation}
M_\rho(x) = \prod_{i=1}^N \left( x_i - \rho^{\frac{s^+}{s_i}}, x_i + \rho^{\frac{s^+}{s_i}} \right),
\end{equation}
where \(s^+ = \max\{s_i : i = 1, 2, \dots, N\}\). For simplicity, we will work primarily with \(M_\rho = M_\rho(0)\).

Next, we define the function space $A(M_\rho, \mathbb{R})$, which is relevant for weak solutions to (1.1), as
\[
A(M_\rho, \mathbb{R}) = \left\{ v : \mathbb{R}^N \to \mathbb{R} : v|_{M_\rho} \in L^p(M_\rho), \quad 
\sum_{i=1}^N \int_{M_\rho} \int_{\mathbb{R}} 
H_i(x, h, |v(x) - v(x + h e_i)|) \frac{dh \, dx}{|h|} < \infty \right\},
\]
where
\begin{equation}
H_i(x, h, \tau) = s_i (1 - s_i) \frac{\tau^p}{|h|^{s_i p}} 
+ t_i (1 - t_i) a(x, x + h e_i) \frac{\tau^q}{|h|^{t_i q}}, 
\quad x \in \mathbb{R}^N, \ h \in \mathbb{R}, \ \tau \geq 0.
\end{equation}

A key tool in our analysis is the nonlocal anisotropic tail, which captures the behavior of solutions at infinity. 
Define the space
\[
L^{q-1}_{\vec{s}, p}(\mathbb{R}) 
= \left\{ v : \mathbb{R}^N \to \mathbb{R} : 
\underset{x \in M_\rho}{\mathrm{ess\,sup}} 
\sum_{i=1}^N \int_{\mathbb{R}} 
\frac{|v(x + h e_i)|^{q-1}}{(1 + |h + x_i|)^{1 + s_i p}} \, dh < \infty \right\}.
\]
Let $m_i \in \{s_i, t_i\}$ for $i = 1, \dots, N$ and $l \in \{p, q\}$. 
Since
\[
\frac{1 + |h + x_i|}{|h + x_i|} \leq 1 + \frac{1}{\rho^{\frac{s^+}{s_i}}}
\quad \text{and} \quad
\frac{|v(x + h e_i)|^{l-1}}{(1 + |h + x_i|)^{1 + m_i l}} 
\leq \frac{|v(x + h e_i)|^{q-1}}{(1 + |h + x_i|)^{1 + s_i p}},
\]
for all $x \in \mathbb{R}^N$ and $|h + x_i| \geq \rho^{\frac{s^+}{s_i}}$, we introduce the following tail quantity:
\begin{equation}\label{eq:Ti_def}
T_i(v,a,b)
:= \underset{x \in M_a}{\mathrm{ess\,sup}} 
\int_{|h + x_i| \geq b^{\frac{s^+}{s_i}}}
\left(
\frac{|v(x + h e_i)|^{p-1}}{|h + x_i|^{1 + s_i p}}
+ \|a\|_{\infty}\frac{|v(x + h e_i)|^{q-1}}{|h + x_i|^{1 + t_i q}}
\right) dh,
\end{equation}
for all $a,b>0$ and $i=1,\dots,N$.
In particular, for $a=b=\rho$, we recover $T_i(v,\rho,\rho)<\infty$ whenever 
$v \in L^{q-1}_{\vec{s}, p}(\mathbb{R})$. This quantity is referred to as the nonlocal anisotropic tail associated with the weak solutions of (1.1). 

It is worth emphasizing that most of the earlier studies in this field focused on the symmetric nonlocal operator, namely
\[
\int_{\mathbb{R}^N}\int_{\mathbb{R}^N}\frac{|u(x)-u(y)|^p}{|x-y|^{N+sp}}\,dx\,dy,
\]
which enjoys full symmetry under the interchange of \(x\) and \(y\). This symmetry simplifies many aspects of the analysis and the derivation of a priori estimates. 

In contrast, the anisotropic expression considered in our work,
\[
\sum_{i=1}^{N}\int_{\mathbb{R}^N}\int_{\mathbb{R}}\frac{|u(x)-u(x+he_{i})|^p}{|h|^{1+s_{i}p}}\,dh\,dx,
\]
lacks this symmetry because it only captures differences along the coordinate directions. This asymmetry introduces significant technical challenges in controlling the nonlocal interactions in various directions. To overcome these difficulties, we imposed the additional condition
\begin{equation}
    \underset{x\in M_{\rho}}{\sup}\sum_{i=1}^{N}\int_{|x_{i}+h|<\rho^{\frac{s^{+}}{s_i}}}|u(x+he_{i})|^{q}\,dh<\infty,
\end{equation}
which provides the necessary integrability framework to handle these directional discrepancies and ensure that weak solutions are locally bounded. It is important to note that the non-symmetric, anisotropic case has not been studied as systematically as the symmetric case, thereby justifying the need for this extra hypothesis.

Now, we state our main Theorem
\begin{thm}
    Let the assumptions (1.3), (1.5), and (1.9) hold, and assume further that 
    \begin{equation}\max\left\{ \frac{1}{1 - t^+}, 1 + s^+ p \right\}\leq q\leq p^{*}_{\bar{s}}.\end{equation}
Then, every weak solution \(u \in A(\mathbb{R}^{N}, \mathbb{R}) \cap L^{q-1}_{\vec{s}, p}(\mathbb{R})\) to (1.1) is locally bounded.
\end{thm}
\begin{Rem}
Condition (1.9) is an \emph{anisotropic integrability condition} that compensates for the lack of symmetry in the kernel. Because the operator only sees differences along coordinate lines, we must ensure that the $q$-phase contributions do not become unbounded in any direction. 

Condition (1.10) serves two purposes. First, it ensures that the higher-order term $L_q^t$ does not dominate the regularizing effect of the lower-order term $L_p^s$, which is crucial for the Sobolev embedding to hold. Second, it enforces the subcriticality $q \le p_{\bar{s}}^*$, preventing the nonlinearity from exceeding the critical exponent. In the isotropic case $s_i = s$, $t_i = t$, condition (1.9) becomes redundant for sufficiently integrable solutions, and (1.10) reduces to the familiar requirement $q \le p_s^*$.

The main novelty of this work lies in establishing local boundedness for nonlocal double-phase
equations that are both anisotropic, in the sense that different fractional orders $s_i$ and $t_i$
act along different directions, and direction-dependent, as the operators are non-symmetric and
capture jumps only along coordinate axes. This combination, motivated by anisotropic jump
processes such as jump--CIR models \cite{fri2}, introduces substantial analytical challenges.
To overcome these challenges, we introduce new tools: anisotropic rectangles $M_\rho$ scaled by $s_i$, anisotropic tail spaces $L^{q-1}_{\vec{s},p}(\mathbb{R})$, and a directionally adapted Caccioppoli estimate (Lemma 4.1), which are essential to handle the lack of symmetry and directional dependence.
\end{Rem}
To prove the local boundedness result, we follow the approaches used in \cite{byu2,cha1,cha3}, adapting techniques to address the double-phase and anisotropic nature of the problem. The primary challenge arises from the interplay between the differing fractional exponents \( s_i \) and \( t_i \) in each direction, which governs the growth and regularity of the operator. This double-phase structure introduces a variable nonlinearity that complicates the analysis compared to homogeneous or isotropic settings. We begin by establishing energy estimates to control the fractional derivatives of the solution, leveraging the framework of anisotropic fractional Sobolev spaces tailored to the exponents \( s_i \) and \( t_i \). A key step involves proving some anisotropic embedding and technical estimates, which allows us to apply the Moser iteration technique. Through this iterative process, we incrementally improve the integrability of the solution, ultimately deriving its local boundedness.

However, the combined effects of anisotropy and double-phase behavior create significant technical hurdles. These include reconciling the distinct scaling properties of \( s_i \) and \( t_i \), as well as harmonizing their directional fractional interactions. While we successfully prove local boundedness using Caccioppoli-type estimates and carefully calibrated scaling arguments, the question of local continuity remains unresolved and constitutes an open problem. This limitation underscores the need for new tools or refined regularity theories capable of bridging the gap between boundedness and continuity in such mixed-exponent settings.

The article is structured as follows: Section 2 presents auxiliary results, including embedding theorems and technical estimates. In Section 3, we investigate the existence of solutions to (1.1). Finally, in Section 4, we establish a Caccioppoli-type estimate and prove Theorem 1.1.

\section{Auxiliary results}
This section aims to provide more or less technical results needed later. In particular, we establish Sobolev-type embeddings tailored to our anisotropic setting. We begin by recalling the following fractional Sobolev--Poincar\'e inequality from [\cite{cha3}, Theorem 1.1].  
From now on, we set $p^{*}:=p^{*}_{\bar{s}}$, unless stated otherwise.
\begin{thm}
    Let $s_{i}\in(0,1)$ and $p>1$ satisfying (1.5) for all $i=1,\ldots,N$. Then there exists a constant $C=C(N,p^{*},p^{*}-p)>0$ such that for all $u\in D^{\vec{s},p}(\mathbb{R}^{N},\mathbb{R})$, we have
    \begin{equation*}
        \|u\|_{L^{p^{*}}(\mathbb{R}^{N})}\leq C\|u\|_{D^{\vec{s},p}(\mathbb{R}^{N},\mathbb{R})}.
    \end{equation*}
\end{thm}

Next, we employ a suitable cutoff function to establish a localized variant of the previously derived Sobolev-type inequality. This refinement allows us to control the function's behavior within a specific subdomain while preserving the essential inequality structure.
\begin{lem}
     Let $s_{i}\in(0,1)$ and $p>1$ satisfying (1.5) for all $i=1,..,N$. Let $M_{\rho}=M_{\rho}(0)\subset \Omega$ be a fixed rectangle with $0<\rho\leq 1$. Then, for any $\frac{\rho}{2}\leq \sigma<\sigma'<\rho$, and $u\in D^{\vec{s},p}(M_{\sigma'},\mathbb{R})$, we get
    \begin{equation*}
        \|u\|_{L^{p^{*}}(M_{\sigma})}\leq C\sum_{i=1}^{N}\biggl\{  s_{i}(1-s_{i})\int_{M_{\sigma'}}\int_{|x_{i}+h|<(\sigma')^{\frac{s^{+}}{s_{i}}}}\frac{|u(x)-u(x+he_{i})|^{p}}{|h|^{1+s_{i}p}}~dhdx + \rho^{-s^{+}p}\|u\|^{p}_{L^{p}(M_{\sigma'})}       \biggr\}^{\frac{1
        }{p}}.
    \end{equation*}
\end{lem}
\begin{proof}
By using Theorem 2.1, and choosing a cut-off function $\tau\in C^{\infty}_{0}(M_{\sigma'})$ satisfying $0\leq\tau\leq1$, $\tau\equiv1$ in $M_{\sigma}$ and $|\frac{\partial\tau}{\partial x_{i}}|\leq\frac{C}{(\sigma'-\sigma)^{\frac{s^{+}}{s_{i}}}}$ for all $i=1,..,N$, we get
\begin{equation}
\begin{split}
     \|u\tau\|_{L^{p^{*}}(\mathbb{R}^{N})}&\leq C\sum_{i=1}^{N}\biggl\{s_{i}(1-s_{i})\int_{\mathbb{R}^{N}}\int_{\mathbb{R}}\frac{|u(x)\tau(x)-u(x+he_{i})\tau(x+he_{i})|^{p}}{|h|^{1+s_{i}p}}~dhdx   \biggr\}^{\frac{1}{p}}\\
     \leq &C\sum_{i=1}^{N}\biggl\{\int_{M_{\sigma'}}\int_{|x_{i}+h|<(\sigma')^{\frac{s^{+}}{s_{i}}}}\frac{|u(x)\tau(x)-u(x+he_{i})\tau(x+he_{i})|^{p}}{|h|^{1+s_{i}p}}~dhdx\\
     &+ \int_{M_{\sigma'}}\int_{|x_{i}+h|\geq(\sigma')^{\frac{s^{+}}{s_{i}}}}\frac{|u(x)\tau(x)|^{p}}{|h|^{1+s_{i}p}}~dhdx
     \biggr\}^{\frac{1}{p}}=C\sum_{i=1}^{N}\{l_{1,i}+l_{2,i}\}^{\frac{1}{p}}.
\end{split}
\end{equation}
By using the following inequality
\begin{equation*}
    |ax-by|\leq \frac{|a-b|}{2}|x+y|+\frac{|a+b|}{2}|x-y|,
\end{equation*}
we get that
\begin{equation}
\begin{split}
    l_{1,i}\leq&\frac{1}{2^{p}}\biggl\{  \int_{M_{\sigma'}}\int_{|x_{i}+h|<(\sigma')^{\frac{s^{+}}{s_{i}}}}2^{p-1}\frac{|u(x)-u(x+he_{i})|^{p}}{|h|^{1+s_{i}p}}|\tau(x)+\tau(x+he_{i})|^{p}  ~dhdx\\
    &+\int_{M_{\sigma'}}\int_{|x_{i}+h|<(\sigma')^{\frac{s^{+}}{s_{i}}}}2^{p-1}\frac{|u(x)+u(x+he_{i})|^{p}}{|h|^{1+s_{i}p}}|\tau(x)-\tau(x+he_{i})|^{p}  ~dhdx
\biggr\}\\
&=\frac{1}{2^{p}}(J_{1,i}+J_{2,i}).
\end{split}
\end{equation}
Thereafter, since $|\tau(x)-\tau(x+he_{i})|^{p}\leq 2^{p}$, we get
\begin{equation}
    J_{1,i}\leq 2^{p}\int_{M_{\sigma'}}\int_{|x_{i}+h|<(\sigma')^{\frac{s^{+}}{s_{i}}}}\frac{|u(x)-u(x+he_{i})|^{p}}{|h|^{1+s_{i}p}}~dhdx,
\end{equation}
and
\begin{equation}
\begin{split}
   J_{2,i}\leq& \int_{M_{\sigma'}}\int_{|x_{i}+h|<(\sigma')^{\frac{s^{+}}{s_{i}}}}\frac{|u(x)|^{p}|\tau(x)-\tau(x+he_{i})|^{p}}{|h|^{1+s_{i}p}}~dhdx\\
   &+2^{p-1}\int_{M_{\sigma'}}\int_{|x_{i}+h|<(\sigma')^{\frac{s^{+}}{s_{i}}}}\frac{|u(x+he_{i})|^{p}|\tau(x)-\tau(x+he_{i})|^{p}}{|h|^{1+s_{i}p}}~dhdx\\
   \leq&2^{p-1}\int_{M_{\sigma'}}|u(x)|^{p}\int_{|x_{i}+h|<(\sigma')^{\frac{s^{+}}{s_{i}}}} \frac{|\tau(x)-\tau(x+he_{i})|^{p}}{|h|^{1+s_{i}p}}~dhdx\\
   &+4^{p-1}\biggl\{\int_{M_{\sigma'}}\int_{|x_{i}+h|<(\sigma')^{\frac{s^{+}}{s_{i}}}}\frac{|u(x)-u(x+he_{i})|^{p}}{|h|^{1+s_{i}p}}|\tau(x)-\tau(x+he_{i})|^{p}~dhdx\\
   &+\int_{M_{\sigma'}}\int_{|x_{i}+h|<(\sigma')^{\frac{s^{+}}{s_{i}}}}\frac{|u(x)|^{p}|\tau(x)-\tau(x+he_{i})|^{p}}{|h|^{1+s_{i}p}}~dhdx\biggr\}\\
   \leq& C\biggl\{  \int_{M_{\sigma'}}|u(x)|^{p}dx\underset{x\in M_{\sigma'}}{\sup}\int_{|x_{i}+h|<(\sigma')^{\frac{s^{+}}{s_{i}}}}\frac{|\tau(x)-\tau(x+he_{i})|^{p}}{|h|^{1+s_{i}p}}~dh\\
   &+\int_{M_{\sigma'}}\int_{|x_{i}+h|<(\sigma')^{\frac{s^{+}}{s_{i}}}}\frac{|u(x)-u(x+he_{i})|^{p}}{|h|^{1+s_{i}p}}~dhdx\biggr\}\\
   \leq& C\biggl\{ \int_{M_{\sigma'}}\int_{|x_{i}+h|<(\sigma')^{\frac{s^{+}}{s_{i}}}}\frac{|u(x)-u(x+he_{i})|^{p}}{|h|^{1+s_{i}p}}~dhdx+\rho^{-s^{+}p} \|u\|^{p}_{L^{p}(M_{\sigma'})}   \biggr\},
\end{split}
\end{equation}
where we used the fact that $|\tau(x)-\tau(x+he_{i})|^{p}\leq\left \|\frac{\partial\tau}{\partial x_{i}}\right\|^{p}_{L^{\infty}}|h|^{p}$, such that
\begin{equation}
\begin{split}
    \underset{x\in M_{\sigma'}}{\sup} \int_{|x_{i}+h|<(\sigma')^{\frac{s^{+}}{s_{i}}}} &\frac{|\tau(x)-\tau(x+he_{i})|^{p}}{|h|^{1+s_{i}p}}~dh \\
    &\leq \left\|\frac{\partial\tau}{\partial x_{i}}\right\|^{p}_{L^{\infty}} \underset{x\in M_{\sigma'}}{\sup} \int_{|x_{i}+h|<(\sigma')^{\frac{s^{+}}{s_{i}}}} |h|^{-1+p(1-s_{i})}~dh \\
    &\leq \frac{C}{(\sigma' - \sigma)^{\frac{s^+}{s_i} p}} (\sigma')^{\frac{s^+}{s_i} p (1 - s_i)} \leq C \rho^{-s^+ p}.
\end{split}
\end{equation}
Next, for $l_{2,i}$, we have
\begin{equation}
    \begin{split}
        l_{2,i}\leq \int_{M_{\sigma'}}|u(x)|^{p}\int_{|x_{i}+h|\geq (\sigma')^{\frac{s^{+}}{s_{i}}}}\frac{1}{|h|^{1+s_{i}p}}~dhdx\leq C\rho^{-s^{+}p}\|u\|^{p}_{L^{p}(M_{\sigma'}).}
    \end{split}
\end{equation}
Therefore, by combining (2.2)-(2.5) into (2.1), we arrive at the desired estimates.
\end{proof}
Next, by defining the average of \( u \) over \( M_{\rho} \) as  
\[
    (u)_{M_{\rho}} = \frac{1}{|M_{\rho}|} \int_{M_{\rho}} u(x) \, dx,
\]
we establish the following local Poincaré inequality for functions in the fractional anisotropic Sobolev spaces.
\begin{lem}
    Let the assumptions and notations of Lemma 2.2 hold. For any function $u \in D^{\vec{s},p}(M_{\sigma'},\mathbb{R})$, we have 
    \begin{equation*}
        \left(\fint_{M_{\sigma}} \left| u - (u)_{M_{\sigma}} \right|^{p^{*}} \,dx\right)^{\frac{p}{p^{*}}} 
        \leq C \rho^{s^{+}p} \sum_{i=1}^{N} \fint_{M_{\sigma'}} \int_{|x_{i}+h|<(\sigma')^{\frac{s^{+}}{s_{i}}}} 
        \frac{|u(x) - u(x + h e_{i})|^{p}}{|h|^{1+s_{i}p}} \,dh \,dx.
    \end{equation*}
\end{lem}
\begin{proof}
By Lemma 2.2, we have
\begin{equation}
\begin{split}
    \|u-(u)_{M_{\sigma}}\|^{p}_{L^{p^{*}}(M_{\sigma})}&\leq C\sum_{i=1}^{N}\int_{M_{\sigma'}}\int_{|x_{i}+h|<(\sigma')^{\frac{s^{+}}{s_{i}}}} \frac{|u(x)-u(x+he_{i})|^{p}}{|h|^{1+s_{i}p}}~dhdx\\
    &~~~~~+ C\rho^{-s^{+}p}\|u-(u)_{M_{\sigma}}\|^{p}_{L^{p}(M_{\sigma'})}.
\end{split}
\end{equation}
Next, we are going to estimate the second term on the right-hand side of (2.7), such that
\begin{equation}
\begin{split}
        \int_{M_{\sigma'}}&|u-(u)_{M_{\sigma}}|^{p}~dx= \int_{M_{\sigma'}}\left|\frac{1}{|M_{\rho}|}\int_{M_{\sigma}}(u(x)-u(y))~dy   \right|^{p}~dx\\
        &\leq C\frac{1}{|M_{\rho}|}\int_{M_{\sigma'}}\int_{M_{\sigma'}}|u(x)-u(y)|^{p}dydx,
\end{split}
\end{equation}
where we used Jensen on the last inequality. Thereafter, we let $l=(l_{0}(x,y),..,l_{N}(x,y))$ be a polygonal chain connecting $x$ and $y$ with
\begin{equation*}
    l_{i}(x,y)=(l^{i}_{1},..,l^{i}_{N}),~\text{where}~\begin{cases}
        l_{j}^{i}=y_{j},~~\text{if}~j\leq i,\\
        l_{j}^{i}=x_{j},~~\text{if}~j> i.
    \end{cases}
\end{equation*}
Then, 
\begin{equation}
\begin{split}
        \int_{M_{\sigma'}}&|u-(u)_{M_{\sigma}}|^{p}~dx\leq \frac{C}{|M_{\rho}|}\sum_{i=1}^{N}\int_{M_{\sigma'}}\int_{M_{\sigma'}}|u(l_{i-1}(x,y))-u(l_{i}(x,y))|^{p}~dydx= \frac{C}{|M_{\rho}|}\sum_{i=1}^{N} I_{i}.
\end{split}
\end{equation}
Next, we fix $i\in \{1,..,N\}$ and set $\omega=l_{i-1}(x,y)=(y_{1},..,y_{i-1},x_{i},..,x_{N})$. Let $z=x+y-\omega=(x_{1},..,x_{i-1},y_{i},..,y_{N})$. Then, $l_{i}(x,y)=\omega+e_{i}(z_{i}-\omega_{i})=(y_{1},..,y_{i},x_{i+1},..,x_{N})$. By Fubini's Theorem, we have
\begin{equation}
\begin{split}
  I_{i}=&\left(\prod_{j\neq i}\int_{-(\sigma')^{\frac{s^{+}}{s_{j}}}}^{(\sigma')^{\frac{s^{+}}{s_{j}}}}~dz_{j} \right)\int_{-(\sigma')^{\frac{s^{+}}{s_{1}}}}^{(\sigma')^{\frac{s^{+}}{s_{1}}}}...\int_{-(\sigma')^{\frac{s^{+}}{s_{N}}}}^{(\sigma')^{\frac{s^{+}}{s_{N}}}}\int_{-(\sigma')^{\frac{s^{+}}{s_{i}}}}^{(\sigma')^{\frac{s^{+}}{s_{i}}}}|u(\omega)-u(\omega+e_{i}(z_{i}-\omega_{i}))|^{p}~dz_{i}d\omega\\
  &\leq C \rho^{\sum_{j\neq i}\frac{s^{+}}{s_{j}}}\int_{M_{\sigma'}}\int_{-(\sigma')^{\frac{s^{+}}{s_{i}}}}^{(\sigma')^{\frac{s^{+}}{s_{i}}}}|u(\omega)-u(\omega+e_{i}(z_{i}-\omega_{i}))|^{p}~dz_{i}d\omega.
        \end{split}
\end{equation}
By letting $h=z_{i}-\omega_{i}$, (2.10) becomes
\begin{equation}
\begin{split}
  I_{i}&\leq C \rho^{s^{+}\beta}\rho^{s^{+}p}\int_{M_{\sigma'}}\int_{|x_{i}+h|<(\sigma')^{\frac{s^{+}}{s_{i}}}}\frac{|u(x)-u(x+he_{i})|^{p}}{|h|^{1+s_{i}p}}~dhdx,
   \end{split}
\end{equation}
where $\beta=\sum_{i=1}^{N}\frac{1}{s_{i}}$. Therefore, by combining (2.11) into (2.9), we arrive at
\begin{equation}
    \int_{M_{\sigma'}}|u-(u)_{M_{\sigma}}|^{p}~dx\leq C\rho^{s^{+}p}\sum_{i=1}^{N}\int_{M_{\sigma'}}\int_{|x_{i}+h|<(\sigma')^{\frac{s^{+}}{s_{i}}}}\frac{|u(x)-u(x+he_{i})|^{p}}{|h|^{1+s_{i}p}}~dhdx.
\end{equation}
Thus, by putting (2.12) into (2.7), and the fact that $|M_{\rho}|=2^{N}\rho^{p^{+}\beta}$ and the definition of $p^{*}$, we get the desired result.
\end{proof}
The following lemma is a deduction of the previous lemma. It will be used to prove Theorem 1.1.
\begin{lem}
Assume that the exponents \( \{s_i\}_{i=1}^{N} \), \( \{t_i\}_{i=1}^{N} \), \( p \), and \( q \) satisfy conditions (1.3) and (1.5), and adopt the notations introduced in Lemma 2.2. Additionally, suppose that \( q \leq p^*_{s} \). For any function \(u \in D^{\vec{s},p}(M_{\sigma'},\mathbb{R})\), we have  
\begin{equation*}
    \begin{split}
        \fint_{M_{\sigma}} |u|^{p} + \|a\|_{\infty} |u|^{q} \,dx 
        &\leq C \biggl\{ \left(\frac{| \operatorname{supp} u |}{|M_{\sigma'}|}  \right)^{\frac{\bar{s}p}{N}}  
        \sum_{i=1}^{N} \fint_{M_{\sigma'}} \int_{|x_{i}+h|<(\sigma')^{\frac{s^{+}}{s_{i}}}} 
        \frac{|u(x)-u(x+he_{i})|^{p}}{|h|^{1+s_{i}p}} \,dh\,dx \\
        &\quad + \|a\|_{\infty} \rho^{s^{+}p} \left( \sum_{i=1}^{N} \fint_{M_{\sigma'}} \int_{|x_{i}+h|<(\sigma')^{\frac{s^{+}}{s_{i}}}} 
        \frac{|u(x)-u(x+he_{i})|^{p}}{|h|^{1+s_{i}p}} \,dh\,dx  \right)^{\frac{q}{p}} \\
        &\quad + \left(\frac{| \operatorname{supp} u |}{|M_{\sigma'}|} \right)^{p-1} \fint_{M_{\sigma}} |u|^{p} 
        + \|a\|_{\infty}|u|^{q} \,dx
        \biggr\}.
    \end{split}
\end{equation*}
\end{lem}
\begin{proof} We follow the lines of [\cite{byu2}, Lemma 2.4].\\
By using Holder's inequality and Lemma 2.3, we have the following estimate
\begin{equation}
\begin{split}
    \fint_{M_{\sigma}}|u|^{q}~dx&\leq C \fint_{M_{\sigma}} |u-(u)_{M_{\sigma}}|^{q}~dx+ C|(u)_{M_{\sigma}}|^{q}\\
    &\leq C\left(\fint_{M_{\sigma}} |u-(u)_{M_{\sigma}}|^{p^{*}}~dx \right)^{\frac{q}{p^{*}}}+ C|(u)_{M_{\sigma}}|^{q}\\
    &\leq C\rho^{s^{+}q}\left(\sum_{i=1}^{N}\fint_{M_{\sigma'}}\int_{|x_{i}+h|<(\sigma')^{\frac{s^{+}}{s_{i}}}}\frac{|u(x)-u(x+he_{i})|^{p}}{|h|^{1+s_{i}p}}~dhdx  \right)^{\frac{q}{p}}+ C|(u)_{M_{\sigma'}}|^{q},
\end{split}
\end{equation}
and, by the same method
\begin{equation}
\begin{split}
    \fint_{M_{\sigma}} |u|^{p}~dx\leq& C\left(\frac{|supp~u|}{|M_{\sigma'}|}  \right)^{\frac{\bar{s}p}{N}}\sum_{i=1}^{N}\fint_{M_{\sigma'}}\int_{|x_{i}+h|<(\sigma')^{\frac{s^{+}}{s_{i}}}}\frac{|u(x)-u(x+he_{i})|^{p}}{|h|^{1+s_{i}p}}~dhdx\\
    &+ C|(u)_{M_{\sigma'}}|^{p}.
\end{split}
\end{equation}
Moreover, we have
\begin{equation}
    \begin{split}
        |(u)_{M_{\sigma'}}|^{p}+\|a\|_{\infty}|(u)_{M_{\sigma'}}|^{q}\leq&\left(\frac{|supp~u|}{|M_{\sigma'}|} \right)^{p-1} \fint_{M_{\sigma'}}|u|^{p}~dx\\
        &+\|a\|_{\infty} \left(\frac{|supp~u|}{|M_{\sigma'}|} \right)^{q-1} \fint_{M_{\sigma'}}|u|^{q}~dx\\
        \leq&\left(\frac{|supp~u|}{|M_{\sigma'}|} \right)^{p-1} \fint_{M_{\sigma'}}|u|^{p}+\|a\|_{\infty}|u|^{q}~dx.
    \end{split}
\end{equation}
Thus, by combining (2.13), (2.14), and (2.15), we get the desired estimate  
\end{proof}
\section{Existence of weak solution}
In this section, we will show the existence of a unique weak solution to the following problem
\begin{equation}
    \begin{cases}
        \mathcal{L}u=0~~~\text{in}~~\Omega,\\
        u=g~~~~~~\text{in}~\mathbb{R}^{N}\setminus \Omega,
    \end{cases}
\end{equation}
where the operator $\mathcal{L}$ is defined in (1.2) and $\Omega$ is a bounded open set in $\mathbb{R}^{N}$ and $g\in A(\mathbb{R}^{N},\mathbb{R})$. Next, by taking
\begin{equation*}
    A_{g}(\Omega)=\{v\in A(\mathbb{R}^{N}, \mathbb{R}):~~u=g~\text{a.e. in}~\mathbb{R}^{N}\setminus \Omega\},
\end{equation*}
we give the following definition.
\begin{defn}
A function $u\in A_{g}(\Omega) $ is said to be a weak solution to the problem (1.1) if 
\begin{equation}
\begin{split}
    \sum^{N}_{i=1}&\int_{\mathbb{R}^{N}}\int_{\mathbb{R}}\bigg\{ s_{i}(1-s_{i})\frac{|u(x)-u(x+he_{i})|^{p-2}}{|h|^{1+s_{i}p}}(u(x)-u(x+he_{i}))\\
    &+t_{i}(1-t_{i})a(x,x+he_{i})\frac{|u(x)-u(x+he_{i})|^{q-2}}{|h|^{1+t_{i}q}}(u(x)-u(x+he_{i}))\biggr\}\\
    &(\eta(x)-\eta(x+he_{i}))~dhdx=0,
\end{split}
\end{equation}
for every  $\eta\in A(\mathbb{R}^{N}, \mathbb{R})$ with $\eta\equiv 0$ a.e. in $\mathbb{R}^{N}\setminus \Omega$.
\end{defn}

Next, by the standard argument used in [[11], Theorem 2.3], we establish the existence, uniqueness, and variational characterization of solutions to problem (1.1). Consider the associated energy functional
\begin{equation}\begin{split}
    \mathcal{E}(v)=&\sum_{i=1}^{N}\int_{\mathbb{R}^{N}}\int_{\mathbb{R}}\frac{s_{i}(1-s_{i})}{p}\frac{|v(x)-v(x+he_{i})|^{p}}{|h|^{1+s_{i}p}}\\
    &+\frac{t_{i}(1-t_{i})}{q}a(x,x+he_{i})\frac{|v(x)-v(x+he_{i})|^{q}}{|h|^{1+t_{i}q}}~dhdx,
\end{split}
\end{equation}
defined for $v\in A(\mathbb{R}^{N},\mathbb{R})$. As a result, we state the following Theorem
\begin{thm}
    Assume that the exponents \(\{s_i\}_{i=1}^{N},~\{t_i\}_{i=1}^{N},~p\), and \(q\) satisfy \((1.3)\) and \((1.5)\), the modulating coefficient $a(.,.)$ satisfies (1.4) and let $g\in A(\mathbb{R}^{N}, \mathbb{R})$. Then, there exists a unique minimizer $u\in A_{g}(\Omega)$ of the functional \(\mathcal{E}\). Moreover, a function $u\in A_{g}(\Omega)$ is a minimizer of \(\mathcal{E}\) if and only if it is a weak solution to problem (3.1) in the sense of Definition 3.1. 
\end{thm}
\begin{proof}
The proof is a straightforward application of the Direct Method of the Calculus of Variations. Let $\{u_{k}\}_{k=1}^{\infty}\subset A_{g}(\Omega)$ be a minimizing sequence. Then we have $\|u_{k}-g\|_{D^{\vec{s},p}(\mathbb{R}^{N},\mathbb{R})}<\infty$, which implies that $\{u_{k}-g\}_{k}$ is uniformly bounded in $L^{p^{*}}(\mathbb{R}^{N})$ by Theorem 2.1. Therefore, by [\cite{cha3}, Theorem 2.1], there exists a subsequence $\{u_{k_{j}}-g\}_{j=1}^{\infty}$ and $v\in L^{p}(M_{\rho})$ where $M_{\rho}\subset\Omega$ such that
\begin{equation*}
    \begin{cases}
        u_{k_{j}}-g \rightharpoonup v~~~~\text{weakly in}~~~~D^{\vec{s},p}(\mathbb{R}^{N},\mathbb{R}),\\
        u_{k_{j}}-g \rightarrow v~~~~\text{strongly in}~~L^{p}(M_{\rho}),~~~~~~~~~~~~~~j\rightarrow\infty\\
        u_{k_{j}}-g \rightarrow v~~~~\text{a.e. in}~~~~~~~~~~M_{\rho}.
    \end{cases}
\end{equation*}
We extend $v$ to $\mathbb{R}^{N}$ by letting $v=0$ on $\mathbb{R}^{N}\setminus M_{\rho}$ and set $u=v+g$. Then, $u_{k_{j}}\rightarrow u$ a.e. in $\mathbb{R}^{N}$. Then, by Fatou's lemma, we have
\begin{equation*}
    \mathcal{E}(u)\leq\underset{j\rightarrow\infty}{\lim\inf}~\mathcal{E}(u_{k_{j}}).
\end{equation*}
This means that $u\in A_{g}(\Omega)$ and it is a minimizer of $\mathcal{E}$. Thereafter, the uniqueness follows from the convexity of the functional.

Now, we prove the equivalence between the minimizer and the weak solution.\\
Step 1: From minimizer to weak solution.
Let \(u \in A_g(\Omega)\) be a minimizer of \(\mathcal{E}\). For any admissible test function \(\eta \in A(\mathbb{R}^N, \mathbb{R})\) with \(\eta \equiv 0\) a.e. in \(\mathbb{R}^N \setminus \Omega\), consider the function
\[
\phi(\varepsilon) = \mathcal{E}(u + \varepsilon \eta), \quad \varepsilon \in \mathbb{R}.
\]
Since \(u\) is a minimizer, we have \(\phi'(0) = 0\). Computing the derivative yields
\begin{equation}
\begin{split}
0 = \phi'(0) =& \sum_{i=1}^N \int_{\mathbb{R}^N} \int_{\mathbb{R}} \bigg[ s_i(1-s_i) \frac{|u(x)-u(x+he_i)|^{p-2}(u(x)-u(x+he_i))}{|h|^{1+s_i p}} \\
& \qquad \qquad \times (\eta(x) - \eta(x+he_i)) \\
& + t_i(1-t_i) a(x,x+he_i) \frac{|u(x)-u(x+he_i)|^{q-2}(u(x)-u(x+he_i))}{|h|^{1+t_i q}} \\
& \qquad \qquad \times (\eta(x) - \eta(x+he_i)) \bigg] \, dh \, dx.
\end{split}
\end{equation}
This is exactly the weak formulation (3.1). Hence \(u\) is a weak solution.\\
Step 2: From weak solution to minimizer.
Conversely, let \(u \in A_g(\Omega)\) be a weak solution. Take any \(v \in A_g(\Omega)\) and set \(\eta = v - u\). Then \(\eta \equiv 0\) a.e. in \(\mathbb{R}^N \setminus \Omega\). By the convexity of the functionals \(t \mapsto |t|^p\) and \(t \mapsto |t|^q\), the energy \(\mathcal{E}\) is convex. Therefore, for any \(\varepsilon \in [0,1]\),
\[
\mathcal{E}(u + \varepsilon \eta) \leq (1-\varepsilon) \mathcal{E}(u) + \varepsilon \mathcal{E}(v),
\]
which implies
\[
\frac{\mathcal{E}(u + \varepsilon \eta) - \mathcal{E}(u)}{\varepsilon} \leq \mathcal{E}(v) - \mathcal{E}(u).
\]
Letting \(\varepsilon \to 0^+\), we obtain
\[
\langle \mathcal{E}'(u), \eta \rangle \leq \mathcal{E}(v) - \mathcal{E}(u),
\]
where \(\langle \mathcal{E}'(u), \eta \rangle\) denotes the Gâteaux derivative of \(\mathcal{E}\) at \(u\) in the direction \(\eta\). Since \(u\) is a weak solution, this derivative vanishes (as shown in (3.4)). Thus,
\[
0 \leq \mathcal{E}(v) - \mathcal{E}(u) \quad \text{for all } v \in A_g(\Omega),
\]
so \(u\) is a minimizer.
\end{proof}
\begin{Rem}\label{rem:nonlocal_dirichlet}
The boundary condition in problem \emph{(3.1)} is understood in the
\emph{nonlocal Dirichlet sense}. Since the operator $\mathcal{L}$ is nonlocal,
the value of $\mathcal{L}u(x)$ for $x \in \Omega$ depends on the values of $u$
in the whole space $\mathbb{R}^N$, and in particular on $\mathbb{R}^N\setminus\Omega$.
For this reason, the condition $u=g$ is imposed almost everywhere on
$\mathbb{R}^N\setminus\Omega$, which is a set of positive measure, rather than
as a trace on $\partial\Omega$.

The weak formulation in Definition~\emph{3.1} is obtained as the
Euler--Lagrange equation associated with the energy functional
$\mathcal{E}$ under variations $\eta \in A(\mathbb{R}^N,\mathbb{R})$ satisfying
$\eta \equiv 0$ almost everywhere in $\mathbb{R}^N\setminus\Omega$.
Within this variational framework, no boundary terms arise and no trace
theory or regularity of $\partial\Omega$ is required.
\end{Rem}

\section{Caccioppoli estimates and local boundedness}
We start with the following lemma which state a nonlocal anisotropic Caccioppoli type estimate. We use (1.7), and define further
\begin{equation}
h_i(x, h, \tau) = s_i (1 - s_i) \frac{\tau^{p-1}}{|h|^{s_i p}} + t_i (1 - t_i) a(x, x + h e_i) \frac{\tau^{q-1}}{|h|^{t_i q}}, \quad x \in \mathbb{R}^N, \, h \in \mathbb{R}, \, \tau \geq 0,
\end{equation}
for all $i=1,..,N$.
\begin{lem}
Assume that the exponents \(\{s_i\}_{i=1}^{N},~\{t_i\}_{i=1}^{N},~p\), and \(q\) satisfy \((1.3)\) and \((1.5)\), and let $u\in A(\mathbb{R}^{N},\mathbb{R})\cap L^{q-1}_{\vec{s},p}(\mathbb{R}^{n})$ be a weak solution to (1.1). Then, for any $M_{\rho}\subset\Omega$ and any $\xi\in C_{0}^{\infty}(M_{\rho})$ with $0\leq\xi\leq1$, the following estimate holds
\begin{equation*}
    \begin{split}
        \sum_{i=1}^{N}&\int_{M_{\rho}}\int_{|x_{i}+h|<\rho^{\frac{s^{+}}{s_{i}}}}H_{i}(x,h,\omega_{\pm}(x)-\omega_{\pm}(x+he_{i}))\max\{\xi(x),\xi(x+he_{i})\}^{q}~dhdx\\
        \leq& \sum_{i=1}^{N}\biggl\{\int_{M_{\rho}}\int_{|x_{i}+h|<\rho^{\frac{s^{+}}{s_{i}}}}H_{i}(x,h,\max\{\omega_{\pm}(x),\omega_{\pm}(x+he_{i})\})|\xi(x)-\xi(x+he_{i})|^{q}~dhdx\\
        &+\left(\int_{M_{\rho}}\omega_{\pm}(x)\xi^{q}(x)~dx  \right)\left( \underset{x\in supp~\xi}{\sup}\int_{|x_{i}+h|\geq\rho^{\frac{s^{+}}{s_{i}}}}h_{i}(x,h,\omega_{\pm}(x+he_{i}))~dh \right)\biggr\},
    \end{split}
\end{equation*}
where $C$ is a positive constant depend on the data only and $\omega_{\pm}:=(u-k)_{\pm}$ with $k\geq0$.
\end{lem}
\begin{proof}
We will only prove the case of $\omega_{+}$, the $\omega_{-}$ case is similar. Thereafter, by taking $\eta=\omega_{+}\xi^{q}$, where $\xi$ is a nonnegative test function in $C_{0}^{\infty}(M_{\rho})$, we get
\begin{equation}
\begin{split}
0\geq& \sum_{i=1}^{N}\biggl\{s_{i}(s_{i}-1)\int_{\mathbb{R}^{N}}\int_{\mathbb{R}}\frac{|u(x)-u(x+he_{i})|^{p-2}}{|h|^{1+s_{i}p}}(u(x)-u(x+he_{i}))\\&~~~~~~~~~(\omega_{+}(x)\xi^{q}(x)-\omega_{+}(x+he_{i})\xi^{q}(x+he_{i}))~dhdx\\
&+t_{i}(t_{i}-1)\int_{\mathbb{R}^{N}}\int_{\mathbb{R}}a(x,x+he_{i})\frac{|u(x)-u(x+he_{i})|^{q-2}}{|h|^{1+t_{i}q}}(u(x)-u(x+he_{i}))\\
&~~~~~~~~~(\omega_{+}(x)\xi^{q}(x)-\omega_{+}(x+he_{i})\xi^{q}(x+he_{i}))~dhdx
\biggr\}\\
=&\sum_{i=1}^{N}\biggl\{\biggl(s_{i}(s_{i}-1)\int_{M_{\rho}}\int_{|x_{i}+h|<\rho^{\frac{s^{+}}{s_{i}}}}\frac{|u(x)-u(x+he_{i})|^{p-2}}{|h|^{1+s_{i}p}}(u(x)-u(x+he_{i}))\\&~~~~~~~~~(\omega_{+}(x)\xi^{q}(x)-\omega_{+}(x+he_{i})\xi^{q}(x+he_{i}))~dhdx\\
&+t_{i}(t_{i}-1)\int_{M_{\rho}}\int_{|x_{i}+h|<\rho^{\frac{s^{+}}{s_{i}}}}a(x,x+he_{i})\frac{|u(x)-u(x+he_{i})|^{q-2}}{|h|^{1+t_{i}q}}(u(x)-u(x+he_{i}))\\
&~~~~~~~~~(\omega_{+}(x)\xi^{q}(x)-\omega_{+}(x+he_{i})\xi^{q}(x+he_{i}))~dhdx\biggr)\\
&+\biggl(s_{i}(s_{i}-1)\int_{M_{\rho}}\int_{|x_{i}+h|\geq\rho^{\frac{s^{+}}{s_{i}}}}\frac{|u(x)-u(x+he_{i})|^{p-2}}{|h|^{1+s_{i}p}}(u(x)-u(x+he_{i}))\\&~~~~~~~~~\omega_{+}(x)\xi^{q}(x)~dhdx\\
&+t_{i}(t_{i}-1)\int_{M_{\rho}}\int_{|x_{i}+h|\geq\rho^{\frac{s^{+}}{s_{i}}}}a(x,x+he_{i})\frac{|u(x)-u(x+he_{i})|^{q-2}}{|h|^{1+s_{i}p}}(u(x)-u(x+he_{i}))\\&~~~~~~~~~\omega_{+}(x)\xi^{q}(x)~dhdx\biggr)\biggl\}\\
=&\sum_{i=1}^{N}\left\{I_{i,1}+ I_{i,2}   \right\}
    \end{split}
\end{equation}
Next, we will assume that $u(x)\geq u(x+he_{i})$ for each $i=1,..,N$. Therefore, we obtain
\begin{equation*}
    \begin{split}
        &|u(x)-u(x+he_{i})|^{l-2}(u(x)-u(x+hei))(\omega_{+}(x)\xi^{q}(x)-\omega_{+}(x+he_{i})\xi^{q}(x+he_{i}))\\
        &=\begin{cases}
(\omega_{+}(x)-\omega_{+}(x+he_{i}))^{l-1}(\omega_{+}(x)\xi^{q}(x)-\omega_{+}(x+he_{i})\xi^{q}(x+he_{i}))~~~~&\text{for}~~u(x)\geq u(x+he_{i})\geq k,\\
(u(x)-u(x+he_{i}))^{l-1}\omega_{+}(x)\xi^{q}(x) ~~~~~~~~~~~~~~~~~~~~~~~~~~~~~~~~~~~~~~~~~~~~~~~~~~~~~~~~~~~&\text{for}~~u(x)>k\geq u(x+he_{i}),\\
0~~~~~~~~~~~~~~~~~~~~~~~~~~~~~~~~~~~~~~~~~~~~~~~~~~~~~~~~~~~~~~~~~~~~~~~~~~~~~~~~~~~~~~~~~~~~~~~~~~~~~~~~~~~~~~~~~~&\text{else},
        \end{cases}\\
        &\geq (\omega_{+}(x)-\omega_{+}(x+he_{i}))^{l-1}(\omega_{+}(x)\xi^{q}(x)-\omega_{+}(x+he_{i})\xi^{q}(x+he_{i})).
    \end{split}
\end{equation*}
Then,
\begin{equation}
    \begin{split}
        I_{i,1}\geq~&~ C\int_{M_{\rho}}\int_{|x_{i}+h|<\rho^{\frac{s^{+}}{s_{i}}}}\frac{(\omega_{+}(x)-\omega_{+}(x+he_{i}))^{p-1}}{|h|^{1+s_{i}p}}\\
        &(\omega_{+}(x)\xi^{q}(x)-\omega_{+}(x+he_{i})\xi^{q}(x+he_{i}))~dhdx\\
        +&C\int_{M_{\rho}}\int_{|x_{i}+h|<\rho^{\frac{s^{+}}{s_{i}}}}a(x,x+he_{i})\frac{(\omega_{+}(x)-\omega_{+}(x+he_{i}))^{q-1}}{|h|^{1+t_{i}q}}\\
        &(\omega_{+}(x)\xi^{q}(x)-\omega_{+}(x+he_{i})\xi^{q}(x+he_{i}))dhdx.      
    \end{split}
\end{equation}
By [\cite{cas2}, Lemma 3.1], we get
\begin{equation}
    \xi^{q}(x)\geq \xi^{q}(x+he_{i})-(1+C_{q}\varepsilon_{i})\varepsilon_{i}^{1-q}|\xi(x)-\xi(x+he_{i})|^{q}-C_{q}\varepsilon_{i}\xi^{1}(x).
\end{equation}
we assume that $\omega_{+}(x)\geq \omega_{+}(x+he_{i})$ and $\xi(x)\leq \xi(x+he_{i})$, for all $i=1,..,N$ (otherwise the estimation bellow is trivial) such that, by taking
\begin{equation*}
    \varepsilon_{i}=\frac{1}{\max\{1,2C_q\}}\frac{\omega_{+}(x)-\omega_{+}(x+he_{i})}{\omega_{+}(x)}\in(0,1),
\end{equation*}
we obtain
\begin{equation*}
    \begin{split}
        \xi^{q}(x)\geq&\xi^{q}(x+he_{i})-\frac{1}{2}\frac{\omega_{+}(x)-\omega_{+}(x+he_{i})}{\omega_{+}(x)}\xi^{q}(x+he_{i})\\
        &-C(q)\frac{(\omega_{+}(x)-\omega_{+}(x+he_{i}))^{1-q}}{\omega_{+}^{1-q}(x)}|\xi(x)-\xi(x+he_{i})|^{q}.
    \end{split}
\end{equation*}
Then, by multiplying the resulting by $(\omega_{+}(x)-\omega_{+}(x+he_{i}))^{l-1}\omega_{x}(x)$, we arrive at
\begin{equation}
   \begin{split}
       (\omega_{+}(x)&-\omega_{+}(x+he_{i}))^{l-1}\xi^{q}(x)\omega_{+}(x)\geq(\omega_{+}(x)-\omega_{+}(x+he_{i}))^{l-1}\omega_{+}(x)\\
       &~\max\{\xi(x),\xi(x+he_{i})\}^{q}-\frac{1}{2}(\omega_{+}(x)-\omega_{+}(x+he_{i}))^{l}\max\{\xi(x),\xi(x+he_{i})\}^{q}\\
       &-C(q)\max\{\omega_{+}(x),\omega_{+}(x+he_{i})\}^{l}|\xi(x)-\xi(x+he_{i})|^{q},
   \end{split} 
\end{equation}
where we assumed that $\xi(x)\leq \xi(x+he_{i})$ and $\omega_{+}(x)\geq\omega_{+}(x+he_{i})$. However, when $0=\omega_{+}(x)\geq \omega_{+}(x+he_{i})\geq 0$, or $\omega_{+}(x)\geq \omega_{+}(x+he_{i})$ and $\xi(x)\geq\xi(x+he_{i})$, the estimate above is trivial. Therefore, we arrive at
\begin{equation}
    \begin{split}
        (&\omega_{+}(x)-\omega_{+}(x+he_{i}))^{l-1}(\xi^{q}(x)\omega_{+}(x)-\xi^{q}(x+he_{i})\omega_{+}(x+he_{i}))\\
        &\geq (\omega_{+}(x)-\omega_{+}(x+he_{i}))^{l-1}\left(\omega_{+}(x)\max\{\xi(x),\xi(x+he_{i})  \}^{q}-\xi^{q}(x+he_{i})\omega_{+}(x+he_{i}) \right)\\
        &-\frac{1}{2}(\omega_{+}(x)-\omega_{+}(x+he_{i}))^{l}\max\{\xi(x),\xi(x+he_{i})\}^{q}-C(q)\max\{\omega_{+}(x),\omega_{+}(x+he_{i})\}^{l}\\
        &~~|\xi(x)-\xi(x+he_{i})|^{q}\\
        &\ge\frac{1}{2}(\omega_{+}(x)-\omega_{+}(x+he_{i}))^{l}\max\{\xi(x),\xi(x+he_{i})\}^{q}-C(q)\max\{\omega(x),\omega(x+he_{i})\}^{l}\\
        &~~|\xi(x)-\xi(x+he_{i})|^{q},
    \end{split}
\end{equation}
whenever $\omega_{+}(x)\geq \omega_{+}(x+he_{i})$ and for all $i=1,..,N$. If $\omega_{+}(x)< \omega_{+}(x+he_{i})$ we exchange the roles of $x$ and $x+he_{i}$ and we obtain (4.6) similarly. As a result, we arrive at
\begin{equation}
    \begin{split}
        I_{i,1}\geq &\frac{1}{2}\biggr\{ s_{i}(1-s_{i})\int_{M_{\rho}}\int_{|x_{i}+h|<\rho^{\frac{s^{+}}{s_{i}}}}\frac{(\omega_{+}(x)-\omega_{+}(x+he_{i}))^{p}\max\{\xi(x),\xi(x+he_{i})\}^{q}}{|h|^{1+s_{i}p}}~dhdx\\
        &+t_{i}(1-t_{i})\int_{M_{\rho}}\int_{|x_{i}+h|<\rho^{\frac{s^{+}}{s_{i}}}}a(x,x+he_{i})\frac{(\omega_{+}(x)-\omega_{+}(x+he_{i}))^{q}}{|h|^{1+t_{i}q}}\\
        &~\max\{\xi(x),\xi(x+he_{i})\}^{q}~dhdx\biggr\}\\
        &-C(q)\int_{M_{\rho}}\int_{|x_{i}+h|<\rho^{\frac{s^{+}}{s_{i}}}}\biggl\{s_{i}(1-s_{i})\frac{\max\{\omega(x),\omega(x+he_{i})\}^{p}}{|h|^{1+s_{i}p}}\\
        &+t_{i}(1-t_{i})a(x,x+he_{i})\frac{\max\{\omega(x),\omega(x+he_{i})\}^{q}}{|h|^{1+t_{i}q}}\biggr\}|\xi(x)-\xi(x+he_{i})|^{q}~dhdx.
    \end{split}
\end{equation}
Next, by using the fact that 
\begin{equation*}
    |u(x)-u(x+he_{i})|^{l-2}(u(x)-u(x+he_{i}))\omega_{+}(x)\geq -\omega_{+}(x+he_{i})^{l-1}\omega_{+}(x),
\end{equation*}
$I_{i,2}$ becomes
\begin{equation}
    \begin{split}
        I_{i,2}\geq& ~-C\int_{M_{\rho}}\int_{|x_{i}+h|\geq\rho^{\frac{s^{+}}{s_{i}}}}\biggr\{\frac{\omega_{+}^{p-1}(x+he_{i})}{|h|^{1+s_{i}p}}+a(x,x+he_{i})\frac{\omega_{+}^{q-1}(x+he_{i})}{|h|^{1+t_{i}q}}\biggl\}\omega_{+}(x)\xi^{q}(x)~dhdx\\
        \geq& -C\left(\int_{M_{\rho}}\omega_{+}(x)\xi^{q}(x)~dx\right)\biggl(\underset{x\in supp~\xi}{\sup}\int_{|x_{i}+h|\geq\rho^{\frac{s^{+}}{s_{i}}}}\frac{\omega_{+}^{p-1}(x+he_{i})}{|h|^{1+s_{i}p}}\\
        &+a(x,x+he_{i})\frac{\omega_{+}^{q-1}(x+he_{i})}{|h|^{1+t_{i}q}}~dh\
\biggr),
    \end{split}
\end{equation}
  for all $i=1,..,N$. 

Hence, by combining (4.7) and (4.8) into (4.2) we get the desired estimate.  
\end{proof}
Now, we are ready to prove the local boundedness Theorem 1.1

\begin{proof}
We follow the idea of the proof of [\cite{byu2}, Theorem 1.1].

Let $M_{\rho}=M_{\rho}(0)$ defined in (1.6) be a fixed rectangle with $\rho\leq1$. For $\frac{\rho}{2}\leq \sigma<\sigma'\leq \rho$, and $k>0$, we set
\begin{equation*}
    A^{+}(k,\rho):=\{x\in M_{\rho}:~u(x)\geq k\}.
\end{equation*}
Next, by denoting
\begin{equation*}
    H(\tau):=\tau^{p}+\|a\|_{\infty}\tau^{q},~~\tau\geq0,
\end{equation*}
and using Lemma 2.4 for $\omega_{+,k}=(u-k)_{+}$ instead of $u$, we get
\begin{equation}
    \begin{split}
      \fint_{M_{\rho}}&H(\omega_{+,k}) ~dx\leq\|a\|_{\infty}\rho^{s^{+}q}\left(\sum_{i=1}^{N}\fint_{M_{\sigma}}\int_{|x_{i}+h|<\sigma^{\frac{s^{+}}{s_{i}}}}\frac{|\omega_{+,k}(x)-\omega_{+,k}(x+he_{i})|^{p}}{|h|^{1+s_{i}p}}~dhdx\right)^{\frac{q}{p}} \\
      &+C\left(\frac{|A^{+}(k,\rho)|}{|M_{\rho}|}  \right)^{\frac{\bar{s}p}{N}}\sum_{i=1}^{N}\fint_{M_{\sigma}}\int_{|x_{i}+h|<\sigma^{\frac{s^{+}}{s_{i}}}}\frac{|\omega_{+,k}(x)-\omega_{+,k}(x+he_{i})|^{p}}{|h|^{1+s_{i}p}}~dhdx\\
      &+ C \left(\frac{|A^{+}(k,\rho)|}{|M_{\rho}|}  \right)^{p-1}\fint_{M_{\sigma}} H(\omega_{+,k})~dx\\
      \leq& C\|a\|_{\infty}\rho^{s^{+}q}\biggl(\sum_{i=1}^{N}\fint_{M_{\sigma'}}\int_{|x_{i}+h|<(\sigma')^{\frac{s^{+}}{s_{i}}}}\frac{|\omega_{+,k}(x)-\omega_{+,k}(x+he_{i})|^{p}}{|h|^{1+s_{i}p}}\\
      &~~~~\max\{\xi(x),\xi(x+he_{i})\}^{q}~dhdx\biggr)^{\frac{q}{p}} \\
      &+C\left(\frac{|A^{+}(k,\rho)|}{|M_{\rho}|}  \right)^{\frac{\bar{s}p}{N}}\sum_{i=1}^{N}\fint_{M_{\sigma'}}\int_{|x_{i}+h|<(\sigma')^{\frac{s^{+}}{s_{i}}}}\frac{|\omega_{+,k}(x)-\omega_{+,k}(x+he_{i})|^{p}}{|h|^{1+s_{i}p}}\\
      &~~~\max\{\xi(x),\xi(x+he_{i})\}^{q}~dhdx\\
      &+ C \left(\frac{|A^{+}(k,\rho)|}{|M_{\rho}|}  \right)^{p-1}\fint_{M_{\sigma'}} H(\omega_{+,k})~dx,
    \end{split}
\end{equation}
where $\xi$ is a cutoff function satisfying $\xi\in C^{\infty}_{0}(M_{\frac{\sigma+\sigma'}{2}})$,~$0\leq\xi\leq1$, $\xi\equiv1$ in $M_{\sigma}$ and $\left|\frac{\partial\xi}{\partial x_{i}}\right|\leq \frac{C}{(\sigma'-\sigma)}$. Next, by fixing $0<l<k$, we have for all $x\in A^{+}(k,\rho)\subset A^{+}(l,\rho)$ that
\begin{equation*}
    (u(x)-l)_{+}=u-l,~~~\text{and}~~\omega_{+,l}:=(u(x)-l)\geq\omega_{+,k},
\end{equation*}
which imply
\begin{equation*}
   |A^{+}(k,\rho)|\leq\int_{A^{+}(k,\rho)}\frac{(u-l)_{+}^{p}}{(k-l)^{p}}~dx\leq \frac{1}{(k-l)^{p}} \int_{A^{+}(k,\rho)}H(\omega_{+,l})~dx,
\end{equation*}
and 
\begin{equation*}
    \int_{M_{\sigma}}\omega_{+,l}~dx\leq\int_{M_{\sigma}}\omega_{+,l}\left(\frac{\omega_{+,l}}{(k-l)}  \right)^{p-1}~dx\leq \frac{1}{(k-l)^{p-1}}\int_{M_{\rho}}H(\omega_{+,l})~dx.
\end{equation*}
Next, by using Lemma 4.1, we get
\begin{equation}
    \begin{split}
\sum_{i=1}^{N}& \int_{M_{\sigma'}}\int_{|x_{i}+h|<(\sigma')^{\frac{s^{+}}{s_{i}}}}\frac{|\omega_{+,k}(x)-\omega_{+,k}(x+he_{i})|^{p}}{|h|^{1+s_{i}p}}\max\{\xi(x),\xi(x+he_{i})\}^{q}~dhdx\\
\leq& C\sum_{i=1}^{N}\biggl\{\int_{M_{\sigma'}}\int_{|x_{i}+h|<(\sigma')^{\frac{s^{+}}{s_{i}}}}\frac{\max\{\omega_{+,k}(x),\omega_{+,k}(x+he_{i})\}^{p}|\xi(x)-\xi(x+he_{i})|^{q}}{|h|^{1+s_{i}p}}~dhdx\\
&+\|a\|_{\infty}\int_{M_{\sigma'}}\int_{|x_{i}+h|<(\sigma')^{\frac{s^{+}}{s_{i}}}}\frac{\max\{\omega_{+,k}(x),\omega_{+,k}(x+he_{i})\}^{q}|\xi(x)-\xi(x+he_{i})|^{q}}{|h|^{1+t_{i}q}}~dhdx\\
&+\left(\int_{M_{\sigma'}}\omega_{+,k}(x)\xi^{q}(x)~dx \right)\biggl(\underset{x\in supp~\xi}{\sup}\int_{|x_{i}+h|\geq (\sigma')^{\frac{s^{+}}{s_{i}}}}\frac{\omega^{p-1}_{+,k}(x+he_{i})}{|h|^{1+s_{i}p}}\\
&+\|a\|_{\infty}\frac{\omega^{q-1}_{+,k}(x+he_{i})}{|h|^{1+t_{i}q}}~dh\biggl)\biggr\}\\
=&C\sum_{i=1}^{N}\left\{I_{i,1}+\|a\|_{\infty}I_{i,2}+I_{i,3}   \right\}.
\end{split}
\end{equation}
Thereafter, we simplify the terms in the right-hand of (4.10). We begin with $I_{i,1}$, such that
\begin{equation}
    \begin{split}
        I_{i,1}\leq& \frac{C}{(\sigma'-\sigma)^{q}}\int_{M_{\sigma'}}\int_{|x_{i}+h|<(\sigma')^{\frac{s^{+}}{s_{i}}}}\frac{\max\{\omega_{+,k}(x),\omega_{+,k}(x+he_{i})\}^{p}}{|h|^{1+s_{i}p-q}}~dhdx\\
        \leq& \frac{C}{(\sigma'-\sigma)^{q}}\int_{M_{\sigma'}}\int_{|x_{i}+h|<(\sigma')^{\frac{s^{+}}{s_{i}}}}\frac{\omega^{p}_{+,k}(x)+\omega^{p}_{+,k}(x+he_{i})}{|h|^{1+s_{i}p-q}}~dhdx\\
        \leq& \frac{C}{(\sigma'-\sigma)^{q}}\biggl\{(\sigma')^{\frac{s^{+}}{s_{i}}(q-s_{i}p)}\int_{M_{\sigma'}}\omega_{+,l}(x)~dx+\int_{M_{\sigma'}}\int_{|x_{i}+h|<(\sigma')^{\frac{s^{+}}{s_{i}}}}\frac{\omega^{p}_{+,k}(x+he_{i})}{|h|^{1+s_{i}p-q}}~dhdx\biggr\},
    \end{split}
\end{equation}
where we used the fact that $|\xi(x)-\xi(x+he_{i})|\leq\left|\frac{\partial\xi}{\partial x_{i}} \right||h|$ for all $i=1,..,N$. Next, we are going to simplify the second integral on the right-hand side of (4.11). Therefore, since $u(x+he_{i})>k>l>0$ for all $-x_{i}-(\sigma')^{\frac{s^{+}}{s_{i}}}<h<-x_{i}+(\sigma')^{\frac{s^{+}}{s_{i}}}$ where $x\in M_{\sigma'}$ (else $\omega^{p}_{+,k}(x+he_{i})=0$), we deduce that $u(x)>k>l>0$ also for $h=0$. As a result, we get
\begin{equation}
    \begin{split}
        \int_{M_{\sigma'}}&\int_{|x_{i}+h|<(\sigma')^{\frac{s^{+}}{s_{i}}}}\frac{\omega^{p}_{+,k}(x+he_{i})}{|h|^{1+s_{i}p-q}}~dhdx\\
        &\leq \frac{C}{(k-l)^{p}} \int_{M_{\sigma'}}\int_{|x_{i}+h|<(\sigma')^{\frac{s^{+}}{s_{i}}}}\omega^{p}_{+,l}(x)\frac{\omega^{p}_{+,k}(x+he_{i})}{|h|^{1+s_{i}p-q}}~dhdx\\
        &\leq\frac{C}{(k-l)^{p}} \int_{M_{\sigma'}}\omega^{p}_{+,l}(x)~dx\left( \underset{x\in M_{\sigma'}}{\sup}   \int_{|x_{i}+h|<(\sigma')^{\frac{s^{+}}{s_{i}}}} \frac{\omega^{p}_{+,k}(x+he_{i})}{|h|^{1+s_{i}p-q}}~dh \right)\\
        &\leq\frac{C}{(k-l)^{p}} (\sigma')^{\frac{s^{+}}{s_{i}}(q-1-s_{i}p)}\int_{M_{\sigma'}}\omega^{p}_{+,l}(x)~dx\left( \underset{x\in M_{\sigma'}}{\sup}   \int_{|x_{i}+h|<(\sigma')^{\frac{s^{+}}{s_{i}}}}\omega^{p}_{+,k}(x+he_{i})~dh \right)~\\
        &\leq\frac{C}{(k-l)^{p}} (\sigma')^{\frac{s^{+}}{s_{i}}(q-1-s_{i}p)}\int_{M_{\sigma'}}\omega^{p}_{+,l}(x)~dx,
    \end{split}
\end{equation}
where we used (1.9) and (1.10). Thereafter, by putting (4.12) into (4.11), we arrive at
\begin{equation}
    I_{i,1}\leq \frac{C}{(\sigma'-\sigma)^{q}}\left\{(\sigma')^{\frac{s^{+}}{s_{i}}(q-s_{i}p)}+\frac{(\sigma')^{\frac{s^{+}}{s_{i}}(q-1-s_{i}p)}}{(k-l)^{p}}    \right\}\int_{M_{\sigma'}}\omega^{p}_{+,l}(x)~dx.
\end{equation}
By the same method, we obtain
\begin{equation}
    I_{i,2}\leq \frac{C}{(\sigma'-\sigma)^{q}}\left\{(\sigma')^{\frac{s^{+}}{s_{i}}(q-t_{i})}+\frac{(\sigma')^{\frac{s^{+}}{s_{i}}(q(1-t_{i})-1))}}{(k-l)^{q}}    \right\}\int_{M_{\sigma'}}\omega^{q}_{+,l}(x)~dx.
\end{equation}
Next, we are going to simplify $I_{i,3}$ such that
\begin{equation}
    \begin{split}
        I_{i,3}=&\left(\int_{M_{\sigma'}}\omega_{+,l}(x)\xi^{q}(x)~dx \right)\biggl(\underset{x\in supp~\xi}{\sup}\int_{|x_{i}+h|\geq(\sigma')^{\frac{s^{+}}{s_{i}}}}\frac{\omega^{p-1}_{+,k}(x+he_{i})}{|h|^{1+s_{i}p}}\\
        &+\|a\|_{\infty}\frac{\omega^{q-1}_{+,k}(x+he_{i})}{|h|^{1+t_{i}q}}~dh\biggr)\\
        \leq& C\left(\frac{\sigma'+\sigma}{\sigma'-\sigma} \right)^{1+t_{i}q}\biggl(\underset{x\in M_{\frac{\sigma'+\sigma}{2}}}{\sup}\int_{|x_{i}+h|\geq(\sigma')^{\frac{s^{+}}{s_{i}}}}\frac{\omega^{p-1}_{+,k}(x+he_{i})}{|x_{i}+h|^{1+s_{i}p}}\\
        &+\|a\|_{\infty}\frac{\omega^{q-1}_{+,k}(x+he_{i})}{|x_{i}+h|^{1+t_{i}q}}~dh\biggr)\int_{M_{\sigma'}}\omega_{+,l}(x)\xi^{q}(x)~dx\\
        =& C \left(\frac{\sigma'+\sigma}{\sigma'-\sigma} \right)^{1+t_{i}q} T_{i}\left(\omega_{+,k},\frac{\sigma'+\sigma}{2},\sigma\right)\int_{M_{\sigma'}}\omega_{+.l}(x)~dx,
    \end{split}
\end{equation}
where we used the fact that
\begin{equation}
    \frac{|x_{i}+h|}{|h|}\leq 1+\frac{|x_{i}|}{|x_{i}+h|-|x|}\leq 1+2\frac{\sigma'+\sigma}{\sigma'-\sigma}\leq 3\frac{\sigma'+\sigma}{\sigma'-\sigma}.
\end{equation}
By putting (4.13)-(4.15) into (4.10), we arrive at
\begin{equation}
    \begin{split}
        \sum_{i=1}^{N}&\fint_{M_{\sigma'}}\int_{|x_{i}+h|<(\sigma')^{\frac{s^{+}}{s_{i}}}}\frac{|\omega_{+,k}(x)-\omega_{+,k}(x+he_{i})|^{p}}{|h|^{1+s_{i}p}}\max\{\xi(x),\xi(x+he_{i})\}^{q}~dhdx\\
        \leq&\frac{C}{(\sigma'-\sigma)^{q}}\sum_{i=1}^{N}\left\{(\sigma')^{\frac{s^{+}}{s_{i}}q(1-t_{i})}+\frac{(\sigma')^{\frac{s^{+}}{s_{i}}(q(1-t_{i})-1)}}{\min\{(k-l)^{p},(k-l)^{q}\}}  \right\}\fint_{M_{\sigma'}}H(\omega_{+,l})~dx\\
        &+C\sum_{i=1}^{N}\left( \frac{\sigma'+\sigma}{\sigma'-\sigma} \right)^{1+t_{i}q}T_{i}\left(\omega_{+,k},\frac{\sigma'+\sigma}{2},\sigma\right)\fint_{M_{\sigma'}}H(\omega_{+,l})~dx.
    \end{split}
\end{equation}
Thereafter, by putting (4.17) into (4.9), we obtain
\begin{equation}
    \begin{split}
       \fint_{M_{\sigma'}}&H(\omega_{+,k})~dx\leq \frac{C\|a\|_{\infty}\rho^{s^{+}p} }{(\sigma'-\sigma)^{\frac{q^{2}}{p}}\min\{(k-l)^{q},(k-l)^{\frac{q^{2}}{p}}\}}\sum_{i=1}^{N}(\sigma')^{\frac{s^{+}q}{s_{i}p}(q(1-t_{i})-1)}\\
       &~~\left( \fint_{M_{\sigma'}}H(\omega_{+,l})~dx  \right)^{\frac{q}{p}}\\
       &+\frac{C}{(k-l)^{\frac{q(p-1)}{p}}}\sum_{i=1}^{N}\left(\frac{\sigma'+\sigma}{\sigma'-\sigma}  \right)^{\frac{q}{p}(1+t_{i}q)}T^{\frac{q}{p}}_{i}\left(\omega_{+,k},\frac{\sigma'+\sigma}{2},\sigma'  \right)\left( \fint_{M_{\sigma'}}H(\omega_{+,l})~dx  \right)^{\frac{q}{p}}\\
       &+\frac{C}{(k-l)^{\frac{p^{2}\bar{s}}{N}+p-1}}\sum_{i=1}^{N}\left(\frac{\sigma'+\sigma}{\sigma'-\sigma}  \right)^{N+t_{i}q}T_{i}\left(\omega_{+,k},\frac{\sigma'+\sigma}{2},\sigma'  \right)\left( \fint_{M_{\sigma'}}H(\omega_{+,l})~dx  \right)^{1+\frac{\bar{s}p}{N}}\\
       &+\frac{C}{(k-l)^{p(p-1)}}\left( \fint_{M_{\sigma'}}H(\omega_{+,l})~dx  \right)^{p}.
    \end{split}
\end{equation}
Now, for $j=0,1,..$ and $k_{0}>1$, we set
\begin{equation*}
    \sigma'_{j}:=\frac{\rho}{2}(1+2^{-j}),~~k_{j}:=2k_{0}(1-2^{-1-j}), ~~\text{and}~~y_{j}:=\int_{A^{+}(k_{j},\sigma_{j})}H(\omega_{+,k_{j}})~dx.
\end{equation*}
Since, $H(u)\in L^{1}(\mathbb{R}^{N})$, we deduce that
\begin{equation*}
    y_{0}=\int_{A^{+}(k_{0},\rho)}H(\omega_{+,k_{0}})~dx\rightarrow0~~\text{as}~k_{0}\rightarrow+\infty.
\end{equation*}
Next, by taking $k_{0}$ large enough such that
\begin{equation}
    y_{j}\leq y_{j-1}\leq...\leq y_{0}\leq1,~~j=1,..,N,
\end{equation}
and since 
\begin{equation*}
   T_{i}(\omega_{+,k_{j}},\sigma_{j},\sigma_{j-1})\leq T_{i}(u,\rho,\rho)<\infty, 
\end{equation*}
where we used (1.8) and the fact that $u\in L^{q-1}_{\bar{s},p}(\mathbb{R}^{N})$, we obtain
\begin{equation*}
\begin{split}
 y_{j+1}\leq& C\biggl\{2^{\frac{2q^{2}}{p}j}y_{j}^{\frac{q}{p}}+2^{\frac{q}{p}(t^{+}q+p)j}y_{j}^{\frac{q}{p}} +2^{(\frac{p^{2}\bar{s}}{N}+p+N+t^{+}q-1))j}y_{j}^{1+\frac{\bar{s}p}{N}}+2^{jp(p-1)}y_{j}^{p} \biggr\}\\
 \leq& C2^{\theta j}y_{j}^{1+\lambda},
 \end{split}
\end{equation*}
where $\theta=\frac{q}{p}(3q+p)+N+2p+q+p^{2}$, and $\lambda=\min\{\frac{q}{p}-1,\frac{\bar{s}p}{N},p-1\}$ for some $C>0$ that depend on the data and $\rho$. Finally, by choosing $k_{0}$ large enough such that
\begin{equation*}
  y_{0}\leq C^{\frac{-1}{\lambda}}2^{\frac{-\theta}{\lambda^{2}}}  
\end{equation*}
holds. Then, for $j\rightarrow\infty$, [\cite{giu}, Lemma 7.1] implies
\begin{equation}
    y_{\infty}=\int_{A^{+}(2k_{0},\frac{\rho}{2})}H((u-2k_{0})_{+})~dx=0,
\end{equation}
which means that $u\leq 2k_{0}$ a.e. in $M_{\frac{\rho}{2}}$.\\

Similarly, by applying the same method we used above to $-u$, we can prove that $u$ is also locally bounded from below.
\end{proof}
\section*{Declarations}
This work was supported by the Carnegie Corporation of New York grant (provided through the AIMS Research and Innovation Center).

\end{document}